 \documentclass[12pt]{amsart}

   \newtheorem{lemma}{Lemma}[section]
   \newtheorem{theorem}[lemma]{Theorem}
   
   \newtheorem{prop}[lemma]{Proposition}
   
   \newtheorem{definition}[lemma]{Definition}

\makeatletter
    \newcommand\figcaption{\def\@captype{figure}\caption}
    \newcommand\tabcaption{\def\@captype{table}\caption}
\makeatother

\title[]{}

\author{Weigu Li}
\address[Weigu Li]
{School of Mathematics\\
 Peking University\\
 Beijing 100871, P.R. China }
\email[W.~Li]{weigu@math.pku.edu.cn}

\author{Kening Lu}
\address[K. Lu]
{Department of Mathematics\\
    Brigham Young University\\
    Provo, Utah 84602, USA } \email[K.~Lu]{klu@math.byu.edu}

\title[Rotation Numbers for Random Dynamical Systems on the Circle]
{Rotation Numbers \\for Random Dynamical Systems on the Circle}

\subjclass{Primary: 60H15; Secondary:
 34C35, 58F11, 58F15, 58F36} \keywords{Rotation number, random maps of circle,
 random differential equations,
\\ This work was partially supported by NSF0200961, NSF0401708, and  NSFC10371083 (K. Lu)
and NSFC10531010 and NNSF10525104(W. Li).  }

\begin{document}

\begin{abstract}
In this paper, we study rotation numbers of random dynamical
systems on the circle. We prove the existence of rotation numbers
and the continuous dependence of rotation numbers on the systems.
As an application, we prove a theorem on analytic conjugacy to a
circle rotation.
\end{abstract}

\maketitle


\section{Introduction}

In this paper, we consider a class of random maps of the circle
arising in the study of dynamical systems when randomness or noise
is taken into account.

\vskip0.1in Let $(\Omega, \mathcal{F}, \mathbb{P})$ be a
probability space and $\theta$ be a measurable
$\mathbb{P}$-measure preserving map on $\Omega$. Let $\psi:
S^1\times \Omega \to S^1$ be an orientation preserving random map
of the circle, where $S^1=\mathbb{R}/2\pi\mathbb{Z}$. This random
map generates a forward random dynamical system
\[
\Psi(n, \cdot, \omega)=
\begin{cases} \psi(\cdot, \theta^{n-1}\omega)\circ \cdots \circ
\psi(\cdot, \omega), & n >0 \\
I, & n=0.
\end{cases}
\]
Let $\varphi(x, \omega): \mathbb{R}\times \Omega \to \mathbb{R}$
be a left of $\psi$,
\[
\psi(e^{ix}, \omega)= e^{i\varphi(x, \omega)},
\]
which satisfies
\begin{itemize}
\item [(1)] $\varphi(x+2\pi,\omega)=\phi(x,\omega)+2\pi;$

\item[(2)] $\varphi(x,\omega)$ is monotonic
 increasing with respect to $x$, i.e.,
 \[\varphi(x,\omega)\ge\varphi(y,\omega),\quad\text{for  }x\ge y.\]
\end{itemize}
We use $\phi(n, \omega)x$ to  denote the random dynamical system
generated by $\varphi$. A random map $\varphi$ is called a
continuous random map if $\varphi(x, \omega)$ is continuous in $x$
for each fixed $\omega\in \Omega$.\footnote{It is enough to assume
that $\varphi$ is continuous in $x$ almost surely in $\omega$ }

We assume that
\begin{itemize}
 \item[(3)]  $\varphi(x,\cdot)\in L^1(\Omega,\mathcal{F},
 \mathbb{P}).$
\end{itemize}
\vskip0.1in

Let $L^1(\Omega,M(S^1))$ denote the set of random maps $\varphi(x,
\omega)$ of the circle  satisfying the above conditions (1), (2),
and (3). We introduce a metric in $L^1(\Omega,M(S^1))$ as follows

\[d(\phi_1, \phi_2)=\int_\Omega\sup_{x\in\mathbb{R}}|\phi_1-\phi_2|\,d\mathbb{P}.
\]
Denote by $L^1(\Omega,H(S^1))$ the subset of
$L^1(\Omega,M(S^1))$ of the continuous random maps. \\

Our main results can be summarized as

\vskip0.1in \noindent {\bf Theorem A.} {\it Let $\varphi \in
L^1(\Omega,M(S^1))$.  Then,
\begin{itemize}
\item[(i)] {\bf Existence of Rotation Number:} There exists a
forward $\theta$-invariant set $\tilde\Omega\in\mathcal{F} $ of
full measure and a $L^1$-function
$\rho(\cdot):\tilde\Omega\to\mathbb{R}$ such that
$$\lim_{n\to\infty}\frac{\phi(n,\omega)x}{2\pi n}=\rho(\omega),\,\,\text{ for all  }x\in\mathbb{R},\,\,\omega\in\tilde\Omega,$$
and $\rho(\theta^n\omega)=\rho(\omega)$ for all $n\in\mathbb{N},$
$\rho(\omega)$ is constant when $\theta$ is   ergodic, where
$\phi(n,\omega)x$ is the random dynamical system generated by
$\varphi$.

\item[(ii)] {\bf Continuous Dependence:}
\[\rho:L^1(\Omega,H(S^1))\to L^1(\Omega,\mathcal{F},
\mathbb{P}):\varphi\mapsto\rho(\omega).\] is continuous

\item[(iii)]{\bf Compact Metric Space:} If, in addition, $\Omega$
is a compact metric space, $\varphi(x,\omega)$ and $\theta$ are
continuous, and $\mathbb{P}$ is the unique $\theta$-invariant
probability measure,
 then
\[\lim_{n\to\infty}\frac{\phi(n,\omega)x}{2\pi n}=\rho,\,\,\text{a real
constant for all }x\in\mathbb{R}\,\text{ and all }\omega\in\Omega.
\]
\item[(iv)]{\bf Random ODE on the Circle:} Analogous results hold
for random ODE's on $S^1$
\[
x'=f(x, \theta^t \omega).
\]

\end{itemize}
}

Theorem A  is an extension of the classical results on rotation
numbers of orientation-preserving homeomorphisms of the circle to
orientation-preserving random maps of the circle and to random
ODE's on the circle. The classical results on rotation numbers can
be found
in \cite{KatokHass}. \\

As an application of this theorem, we consider the case when
$\Omega=\mathbb{R}^{m-1}/2\pi\mathbb{Z}$ is the torus of dimension
$m-1$ and the  dynamical system $\{\theta^{n}\}_{n\in\mathbb{Z}}$
on $\Omega$ is given by
\[\theta^n:\omega\mapsto\omega+2n\pi\alpha,\]
where $\alpha\in\mathbb{R}^{m-1}$ is a given vector. We assume
that
$$\langle\alpha,k\rangle-j\ne0,\,\,\text{for all }k\in\mathbb{Z}^{m-1}\backslash\{0\},\,j\in\mathbb{Z}.$$
Then,  the normalized Lebesgue measure $\mathbb{P}$ is the unique
$\theta$-invariant probability measure and $\theta$ is ergodic
under $\mathbb{P}$. Let $\varphi(x,\omega)$ be a random map of the
circle in $L^1(\Omega,M(S^1))$.   Suppose that
$\varphi(\cdot,\cdot)$ is continuous. Then,  by Theorem A, the
rotation number of $\varphi$ exists and is given by
\[\rho=\lim_{n\to\infty}\frac{\phi(n, \omega)x}{2\pi n}=\lim_{n\to\infty}\frac1{2\pi n}
\varphi(\cdot, {\theta^{n-1}\omega})\circ \varphi(\cdot,
{\theta^{n-2} \omega})\circ\cdots\circ \varphi(x,\omega), \quad
\text{for all } x\in \mathbb{R}, \omega\in \Omega, \] which is
independent of $\omega$
and $x$. \\

 We denote by $U_r$ the
strip region in $\mathbb{C}^m$,
\[U_r:=\{z=(z_1,z_2,\dots,z_m)\in\mathbb{C}^m:|\text{Im}
z_i|<r,\,i=1,2,\dots,m\}.
\] where $r$ is a positive number. For a holomorphic function $p(z)$ bounded in this region, we
define
$$\|p\|_r=\sup_{z\in U_r}|p(z)|.$$ For a vector valued analytic
function $f(z)=(f_1(z),f_2(z),\dots,f_m(z)):U_r\to\mathbb{C}^m$,
we define
$$\|f\|_r=\max_{1\le j\le m}\|f_j\|_r.$$

Consider a perturbation of the circle rotation by $2\pi \rho$,
\[
\varphi(x, \omega)=x+2\pi \rho +p(x, \omega).
\]
We have the following theorem on analytic conjugacy to a circle
rotation.

\vskip0.1in

\noindent{ \bf Theorem B. } {\it Let $p(x,\omega)$ be analytic in
$U_r$ and $2\pi$-periodic in each variable, real on the real axes.
Assume
\begin{itemize}

\item[(1)] $\varphi(x,\omega)=x+2\pi\rho+p(x,\omega)$ has the
rotation number $\rho$ and

\item[(2)] the vector $\mu=(\rho,\alpha)$ is of $(C,\,\nu)$ type,
i.e.,
\begin{equation}\label{cnu}|e^{2\pi i\langle\mu,\,k\rangle}-1|>\frac
C{|k|^\nu}, \quad |k|:=|k_1|+|k_2|+\cdots+|k_m|
\end{equation}
for all nonzero integer vector $k\in\mathbb{Z}^m,$ where $C$ and
$\nu$ are positive constants.
\end{itemize}
Then, there exists $\epsilon>0$ depending only on $C,\nu, r$ and
$m$ such that if $\|p\|_r<\epsilon$, then the random map
$\varphi(\cdot, \omega)$ is analytically conjugate to the circle
rotation by the angle $2\pi\rho$, i.e.,  there exists an
analytical random transformation
\[H(\cdot,\cdot):\mathbb{R}^m/2\pi\mathbb{Z}\to\mathbb{R}/2\pi\mathbb{Z}\]
such that
\[H(x+2\pi\rho, \theta\omega)=\varphi (\cdot, \omega)\circ H(x, \omega).\]
}

\vskip0.1in When $\varphi(x, \omega)$ is independent of $\omega$,
Theorem B was first proved by V. Arnold \cite{Arnold1},
\cite{Arnold2}, see also \cite{KatokHass}. Theorem B is also
related to the classical result that a family of analytic
homeomorphisms
\[
y \mapsto y+\alpha+ p(y), \quad\text{where } \alpha,
y\in\mathbb{T}^n,
\]
for most of $\alpha$, is analytically conjugate to a translation
$y \mapsto y+2\pi\mu$ if $p$ is sufficiently small.

\vskip0.1in

The study of rotation numbers of homeomorphims of the circle goes
back to Poincar\'e and has a rich history.  It plays an important
role in the investigation of qualitative behavior of various
dynamical systems. There is an extensive literature on this
subject. We will not try to give an exhaust list of references,
but mention only some related works here. For the classical
results such as the Poincar\'e classification theorem and the
Denjoy theorem we refer to Katok and Hasselblatt \cite{KatokHass}.
For the application of rotation numbers to the spectral theory of
almost periodic Schr\"odinger operators, see Johnson and Moser
\cite{JohnsonMoser}. Recently, Fabbri, Johnson, and Nunez
\cite{Johnson1}, \cite{Johnson2}, \cite{Johnson3} have obtained
the results on rotation numbers for non-autonomous linear
Hamiltonian Systems. The rotation numbers of asymmetric equations
were  studied by Feng and Zhang \cite{Zhang}. For the rotation
numbers of stochastic ODE's, see the book by L. Arnold and the
references therein. The existence of rotation numbers in Theorem A
is also related to the work by Ruffino \cite{Ruffino} on rotation
numbers of linear processes in $\mathbb{R}^2$. There is a nice
survey article by Franks \cite{Franks} on rotation numbers in
dynamics on surfaces and their applications to the description of
dynamical systems. For diffeomorphisms of the circle, there are
rich results on smooth conjugacy to circle rotations, see the
works of Arnold\cite{Arnold2}, Herman \cite{Herman},
Yoccoz\cite{Yoccoz1}, Katznelson and Ornstein \cite{KO},  Sinai
and Khanin \cite{SK},
and their references therein.\\

We organize this paper as follows. In section 2, we prove the main
results on rotation numbers for random maps of the circle.
Analogous results for random ODE's on the circle  are given in
Section 3. In section 4, we apply our main results to almost
periodic ordinary differential equations and random compositions
of homeomorphisms of the circle. The proof of Theorem B and a more
general results are given in Section 5. The proof of Theorem B is
based on the standard KAM approach.

\vskip0.2in \noindent{\bf Acknowledgement.} We would like to thank
Jiangong You for his comments and suggestion on random
compositions of homeomorphisms of the circle.  The first author
also would like to thank the Department of Mathematic at Brigham
Young University for the support and hospitality during the fall
of 2004 when this paper was written.

\section{\bf Rotation Numbers of Random Maps of the Circle}

In this section, we prove our main results for orientation
preserving random maps of the circle. We formulate properties (1),
(2), and (3) in Theorem A as the following propositions.

\begin{prop} \label{maprot1} (Existence of Rotation Numbers)  Let $\varphi \in
L^1(\Omega,M(S^1))$.  Then, there exists a forward
$\theta$-invariant set $\tilde\Omega\in\mathcal{F} $ of full
measure and a $L^1$-function
$\rho(\omega):\tilde\Omega\to\mathbb{R}$ such that
$$\lim_{n\to\infty}\frac{\phi(n,\omega)x}{2\pi n}=\rho(\omega),\,\,\text{ for all }x\in\mathbb{R},\,\,\omega\in\tilde\Omega,$$
and $\rho(\theta^n\omega)=\rho(\omega)$ for all $n\in\mathbb{N},$
$\rho(\omega)$ is constant when $\theta$ is  ergodic, where
$\phi(n,\omega)x$ is the random dynamical system generated by
$\varphi$.
\end{prop}

\begin{proof}
First we claim that
\begin{equation}\label{li1}
\sup_{x\in\mathbb{R}}(\phi(n,\omega)x-x)\le\inf_{x\in\mathbb{R}}(\phi(n,\omega)x-x)+2\pi.
\end{equation}
Indeed, let $y<x<y+2\pi$, then
\[\phi(n,\omega)y-y-2\pi\le\phi(n,\omega)y-x\le\phi(n,\omega)x-x\le\phi(n,\omega)y+2\pi-x\le\phi(n,\omega)y-y+2\pi,
\] which yields (\ref{li1}).
For any $x<y$, let $m=[y-x]+1$, where $[\,\cdot\,]$ denotes the
integer part of the number, then
$$\phi(n,\omega)x\le\phi(n,\omega)y\le\phi(n,\omega)x+2\pi m.$$
Thus, if the limit $\lim_{n\to\infty}\phi(n,\omega)x/n$ exists, it
is independent of the point $x$. Let
\[\Theta(x,\omega):=(\phi(1,\omega)x,\theta\omega)
\]denote the corresponding skew product map and write
$\phi(1,\omega)x=\varphi(x,\omega)=x+h(x,\omega)$. Then,
$h(x+2\pi,\omega)=h(x,\omega)$. Thus, by the definition of
$\Theta$, we have
\begin{align}\label{eq1}\begin{split}
\phi(s+t,\omega)(0)&=\sum^{s+t-1}_{k=0}h\circ\Theta^k(0,\omega)=\phi(s,\omega)(0)+\sum^{s+t-1}_{k=s}h\circ\Theta^k(0,\omega)\\&
=\phi(s,\omega)(0)+\phi(t,\theta^s\omega)\phi(s,\omega)(0)-\phi(s,\omega)(0).
\end{split}
\end{align}
Hence, by  the claim (\ref{li1}), we have
\begin{equation}\label{li2}\phi(s,\omega)(0)+\phi(t,\theta^s\omega)(0)-2\pi\le\phi(s+t,\omega)(0)\le\phi(s,\omega)(0)+\phi(t,\theta^s\omega)(0)+2\pi.
\end{equation}
 Let
$f_n(\omega)=\phi(n,\omega)(0)+2\pi,g_n(\omega)=-\phi(n,\omega)(0)+2\pi,$
then by inequality (\ref{li2}), we obtain
\[f_{s+t}(\omega)\le f_s(\omega)+f_t(\theta^s\omega),\quad g_{s+t}(\omega)\le g_s(\omega)+g_t(\theta^s\omega).\]
By Kingman's Subadditive Ergodic Theorem, there is a forward
invariant set $\tilde\Omega$ of full measure and measurable
functions
$\rho(\omega),\tilde\rho(\omega):\tilde\Omega\to\mathbb{R}\cup\{-\infty\}$
with
$\rho(\theta\omega)=\rho(\omega),\tilde\rho(\theta\omega)=\tilde\rho(\omega)$
and $\rho^+,\tilde\rho^+\in L^1(\Omega,\mathcal{F},\mathbb{P})$,
where $a^+:=\max(0,a),$ such that
\[\lim_{n\to\infty}\frac{f_n(\omega)}{2\pi n}=\rho(\omega),\,\,\quad
\lim_{n\to\infty}\frac{g_n(\omega)}{2\pi
n}=\tilde\rho(\omega),\,\,\text{ for
}x\in\mathbb{R},\,\,\omega\in\tilde\Omega,\] Obviously,
$\rho(\omega)=-\tilde\rho(\omega),$ which implies that
$\rho(\omega)\in L^1(\Omega,\mathcal{F},\mathbb{P})$. When
$\theta$ is ergodic, $\rho$ is a constant. This completes the
proof of the proposition.
\end{proof}

\begin{prop} \label{maprot2} (Continuous Dependence of Rotation
Numbers)   The mapping
\[\rho:L^1(\Omega,H(S^1))\to L^1(\Omega,\mathcal{F},
\mathbb{P}):\varphi\mapsto\rho(\omega).\] is continuous
\end{prop}

Before we prove this proposition, we review some basic concepts
and results on random dynamical systems, which are taken from
Arnold \cite{Arn98}. Let $(\Omega, \mathcal{F}, \mathbb{P})$ be a
probability space. Let $\mathbb{T}=\mathbb{R},\; \mathbb{R}^+,\;
\mathbb{Z},\;$ or $\mathbb{Z}^+$. $\mathbb{T}$ is endowed with its
Borel $\sigma$-algebra $\mathcal{B}(\mathbb{T})$.

\begin{definition} A family $(\theta^t)_{t\in \mathbb{T}}$ of
mappings from $\Omega$ into itself is called a metric dynamical
system if

\begin{itemize}
\item[(1)] $(\omega, t) \to \theta^t\omega$ is $\mathcal{F}\otimes
\mathcal{B}(\mathbb{T})$ measurable;

\item[(2)] $\theta^0=id_{\Omega}$, the identity on $\Omega$ ,
$\theta^{t+s}=\theta^t\circ \theta^s$ for all $t, s \in
\mathbb{T}$;

\item[(3)] $\theta^t$ preserves the probability measure
$\mathbb{P}$.
\end{itemize}
\end{definition}

\begin{definition} Let $X$ be a matric space. A map
$$
\phi: \mathbb{T} \times \Omega \times X \to X, \quad (t, \omega,
x) \mapsto \phi(t, \omega, x),
$$
is called a random dynamical system (or a cocycle) over a metric
dynamical system $(\Omega, \mathcal{F}, \mathbb{P}, {\theta}^t)$
if
\begin{itemize}

\item[(1)] $\phi$ is
$\mathcal{B}(\mathbb{T})\otimes\mathcal{F}\otimes\mathcal{B}(X)$-measurable;

\item[(2)] The map $\phi(t,\omega):=\phi(t,\omega,\cdot): X \to X$
forms a cocycle over $\theta^t$:
\[
\phi(0, \omega)=Id, \quad \hbox{ for all }\; \omega \in \Omega,
\]
\[
\phi(t+s,\omega)=\phi(t,\theta^{s}\omega)\circ \phi(s,\omega),
\quad \hbox{ for all }\; t, s \in \mathbb{T}, \quad\omega \in
\Omega.
\]

\end{itemize}
\end{definition}

A probability measure $\mu$ on $(X\times\Omega
,\mathcal{B}\otimes\mathcal{F})$ is said to be an invariant
measure for the random dynamical system $\phi$, or
$\phi\,$-invariant, if it satisfies
\begin{itemize}
\item [(i)]  $\Theta(t)\mu=\mu$ for all $t\in\mathbb{T},$

\item[(ii)] $\pi_\Omega\mu=\mathbb{P},$
\end{itemize}
where $\Theta(t)(x, \omega)=(\phi(t, \omega) x, \theta^t\omega)$
and $\pi_\Omega: X\times \Omega \to \Omega$ is the projection on
$\Omega$.

Let
$$\mathcal{P}_{\mathbb{P}}(X\times \Omega):=\{\mu\;:\; \mu \; \text{is a probability measure on } (X\times \Omega,\mathcal{B}\otimes\mathcal{F})\text{ with marginal }
\mathbb{P} \text{ on } (\Omega,\mathcal{F})\}$$ and $$
\mathcal{I}_\mathbb{P}(\phi):=\{\mu\in\mathcal{P}_{\mathbb{P}}(X\times
\Omega):\mu\;\text{is } \phi\,\text{-invariant}\}.
$$ By Theorem 1.5.10 in \cite{Arn98}, if $X$ is a compact metric space and $\phi$
is a continuous random dynamical system on  $X$, then
$\mathcal{I}_\mathbb{P}(\phi)$ is not empty.

Suppose that $X$ is a compact metric space, by Proposition 1.4.3
in \cite{Arn98}, any probability measure
$\mu\in\mathcal{P}_{\mathbb{P}}(X\times \Omega)$ has a
$\mathbb{P}$-a.s. unique factorization
$$\mu(dx,d\omega)=\mu_\omega(dx)\mathbb{P}(d\omega),$$
or equivalently: for all $f\in L^1_\mu(X\times \Omega)$
$$\int_{X\times\Omega}f\,d\mu=\int_\Omega\left(\int_Xf(x,\omega)\mu_\omega(dx)\right)\mathbb{P}(d\omega).$$
Let $\mathcal{C}_b(X)$ denote the Banach space of real-valued
bounded continuous functions on $X$, with sup norm
$\|f\|_b:=\sup_{x\in X}|f(x)|$. We call a function
$f:\Omega\to\mathcal{C}_b(X)$ measurable if $(x,\omega)\mapsto
f(x,\omega)$ is measurable. Define
$$L^1_{\mathbb{P}}(\Omega,\mathcal{C}_b(X)):=\{f:\Omega \to\mathcal{C}_b(X)
\text{ measurable,
}\|f\|:=\int_\Omega\|f(\omega,\cdot)\|_bd\mathbb{P}<\infty\}.$$
\begin{definition}\cite{Arn98} We call the smallest topology in $\mathcal{P}_{\mathbb{P}}(X\times \Omega)$ which makes
$\mu\mapsto \mu(f)$ continuous for each $f\in
L^1_{\mathbb{P}}(\Omega,\mathcal{C}_b(X))$ the weak convergence on
$\mathcal{P}_{\mathbb{P}}(X\times \Omega)$. A sequence $\{\mu_i\}$
converges in the topology to $\mu$ if $\mu_i(f)\to\mu(f)$ for each
$f\in L^1_{\mathbb{P}}(\Omega,\mathcal{C}_b(X))$.
\end{definition}
Let $X$ be a metric space and  $ L^1_{\mathbb{P}}(\Omega,{\mathcal
L i p}_b(X))$ denote the subset of $
L^1_{\mathbb{P}}(\Omega,\mathcal{C}_b(X))$ such that
$$Lip(f):=\left|\sup_{x\ne y}\frac{|f(x,\cdot)-f(y,\cdot)|}{|x-y|}\right|_{L^\infty_{\mathbb{P}}(\Omega)}<\infty.$$
\begin{lemma}\label{mea}
 Let $X$ be a metric space and $\mu,\nu\in \mathcal{P}_{\mathbb{P}}(X\times
 \Omega)$. Then
 $$\mu=\nu\Leftrightarrow\int f\,d\mu=\int fd\nu,\,\,\,\text{for all } f\in L^1_{\mathbb{P}}(\Omega,{\mathcal Lip}_b(X)).$$
\end{lemma}
The proof of this lemma follows the same lines as Lemma 1.5.4.
\cite{Arn98}
\begin{lemma}\label{weak}
 Let $X$ be a compact metric space. Let $\phi_i(t,\omega)x$ be a
sequence of continuous random dynamical systems over the metric
dynamical system  \{$\theta^t\}_{t\in\mathbb{T}}$  on $X$ and
$\mu_i$ be an invariant measure of $\phi_i$. If $\phi_i$ converges
to a continuous random dynamical system $\phi(t,\omega)x$ in the
following sense:
$$\int_\Omega \sup_{x\in X}d(\phi_i(t,\omega)x,\phi(t,\omega)x)\,d\mathbb{P}\to0,\text{ for any given }t\in\mathbb{T}\cap[0,1],$$
then every limit point of $\mu_i$ for $i\to\infty$ in the topology
of weak convergence is an invariant measure for $\phi$.
\end{lemma}
\begin{proof}Assume that $\mu_i\to\mu$ in the topology of weak
convergence. Since $\mathcal{P}_{\mathbb{P}}(X\times \Omega)$ is
compact by Theorem 1.5.10 in \cite{Arn98}, $\mu\in
\mathcal{P}_{\mathbb{P}}(X\times \Omega)$. Next, we show that
$\mu$ is invariant for $\phi$. Denote by $\Theta_i(t)$ and
$\Theta(t)$ the corresponding skew-product of $\phi_i$ and $\phi$
respectively. For fixed $t\in\mathbb{T}\cap[0,1]$ and for any
given $f\in L^1_{\mathbb{P}}(\Omega,{\mathcal L i p}_b(X))$,
\[\lim_{i\to\infty}|(\Theta_i(t)\mu_i(f)-\Theta(t)\mu_i(f)|\le\lim_{i\to\infty}
Lip(f)\int_\Omega \sup_{x\in
X}d(\phi_i(t,\omega)x,\phi(t,\omega)x)\, d\mathbb{P}=0.\]
Therefore,
\[\Theta(t)\mu(f)=\lim_{i\to\infty}
\Theta(t)\mu_i(f)=\lim_{i\to\infty}\Theta_i(t)\mu_i(f)=\lim_{i\to\infty}\mu_i(f)=\mu(f).\]
This implies by Lemma \ref{mea},  $\Theta(t)\mu=\mu$. The proof is
complete.

\end{proof}

\noindent{\bf Proof of Proposition \ref{maprot2}.}
 Denote by
$\rho_\varphi(\omega)$ the rotation number of $\varphi$. First, we
note that the rotation number $\rho$ is monotonic with respect to
$\varphi$, i.e.,  $\rho_\varphi\ge\rho_{\tilde\varphi}$ for
$\varphi\ge\tilde\varphi.$ If the proposition is not true, there
exists a sequence of mappings
$\varphi_i(x,\omega)=x+h_i(x,\omega)\in L^1(\Omega,H(S^1))$  with
$\varphi_i\to\varphi=x+h(x,\omega)$ in $L^1(\Omega,H(S^1))$ such
that
\begin{equation}\label{ne0}\lim_{i\to\infty}\int_\Omega|\rho_{\varphi_i}(\omega)-\rho_\varphi(\omega)|d\mathbb{P}\ne0.\end{equation}
 Let
$$\varphi^+_i=x+h^+_i(x,\omega):=\max\{\varphi,\varphi_i\},\,\,\phi^-_i=x+h^-_i(x,\omega):
=\min\{\varphi,\varphi_i\},$$ then
$\|\varphi^+_i-\varphi^-_i\|\to0.$ Let $\nu^+_i$ and $\nu^-_i$ be
the invariant probability  measures of $\phi^+_i$ and $\phi^-_i$
respectively. Then by the fact that
$\mathcal{P}_{\mathbb{P}}(X\times \Omega)$ is compact and Lemma
\ref{weak}, we can assume that
 they converge to some
invariant measures $\nu^+$ and $\nu^-$ of $\phi$ in the weak
topology respectively. Finally, by using Fubini's Theorem,
(\ref{eq1}), and Birkhoff's Ergodic Theorem, we have
\begin{align*}&\lim_{i\to\infty}\int_\Omega|\rho_{\varphi_i}(\omega)-\rho_\varphi(\omega)|d\mathbb{P}\\
&\le\lim_{i\to\infty}
\int_\Omega(\rho_{\varphi^+_i}(\omega)-\rho_{\varphi^-_i}(\omega))d\mathbb{P}\\
&=\lim_{i\to\infty}
\left(\int_{X\times\Omega}\rho_{\phi^+_i}(\omega)d\nu^+_i-\int_{X\times\Omega}
\rho_{\phi^-_i}(\omega)d\nu^-_i\right)
\\&=\lim_{i\to\infty}\left(\int_{X\times\Omega} h^+_id\nu^+_i-\int_{X\times\Omega} h^-_id\nu^-_i\right)\\&=\lim_{i\to\infty}\left(\int_{X\times\Omega}
 (h^+_i-h)d\nu^+_i-\int_{X\times\Omega}
(h^-_i-h)d\nu^-_i+\int_{X\times\Omega}
h\,d\nu^+_i-\int_{X\times\Omega}
h\,d\nu^-_i\right)\\&\le\lim_{i\to\infty}(\|h-h^+_i\|+\|h-h^-_i\|)+\lim_{i\to\infty}\left(\int_{X\times\Omega}
h\,d\nu^+_i-\int_{X\times\Omega}
h\,d\nu^-_i\right)\\
&=\int_{X\times\Omega} h\,d\nu^+-\int_{X\times\Omega}
h\,d\nu^-\\&=\int_\Omega\rho_\varphi
d\mathbb{P}-\int_\Omega\rho_\varphi d\mathbb{P}=0,
\end{align*}    which contradicts to (\ref{ne0}). This completes
the proof of Proposition \ref{maprot2}.

\hfill $\square$

\vskip0.1in

\begin{prop} \label{maprot3} (Compact Metric Space)  Let $\varphi \in
L^1(\Omega,M(S^1))$. If $\Omega$ is a compact metric space,
$\varphi(x,\omega)$ and $\theta$ are continuous, and $\mathbb{P}$
is the unique $\theta$-invariant probability measure,
 then
\[\lim_{n\to\infty}\frac{\phi(n,\omega)x}{2\pi n}=\rho,\,\,\text{a real
constant for all }x\in\mathbb{R},\,\,\text{all }\omega\in\Omega.
\]
\end{prop}

\begin{proof} Since $\Omega$ is compact, on the one hand, according to
Krylov and Bogolyubov's Theorem, the skew product map
$\Theta(x,\omega):S^1\times\Omega\to S^1\times\Omega$ possesses at
least one invariant probability  measure. On the other hand, by
the assumption that $\mathbb{P}$ is the unique $\theta$-invariant
probability measure, we have that $\theta$ is ergodic and
$\pi_\Omega\mu=\mathbb{P}$ for any invariant probability  measure
$\mu$ of $\Theta$,   where
$\pi_\Omega:S^1\times\Omega\to\Omega,\,\pi_\Omega(x,\omega)=\omega$,
is the projection onto $\Omega$. Therefore, using (\ref{eq1}) for
any invariant probability  measure $\mu$ of $\Theta$, we have
\[\lim_{n\to\infty}\frac{\phi(n,\omega)x}{2\pi n}=\lim_{n\to\infty}\frac1{2n\pi}\sum^{n-1}_{k=0}h\circ\Theta^k(x,\omega)=\rho,\,\mu\text{-}a.s.\]
where $\rho$ is a real number. Thus, by Birkhoff's Ergodic
Theorem, we obtain
\begin{equation}\label{li3}\int_{S^1\times
\Omega} h\,d\mu=\lim_{n\to\infty}\int_{S^1\times
\Omega}\frac1n\sum^{n-1}_{k=0}h\circ\Theta^k\,d\mu =\int_{S^1\times
\Omega}\lim_{n\to\infty}\frac1n\sum^{n-1}_{k=0}h\circ\Theta^k\,d\mu=2\pi\rho.
\end{equation}
To complete the proof of Proposition \ref{maprot3}, we need the
following lemma.
\begin{lemma}\label{mean} Let $B$ be a compact metric space
or a compact Hausdorff separable space and $
\{\Theta^t\}_{t\in\mathbb{T}}$ be a continuous dynamical system on
$B$. Suppose that $G$ is a continuous function on $B$ such that
$$\int_B G\,d\mu=0$$
for all invariant probability  measure $\mu$ of $\Theta$. Then
\[\lim_{n\to\infty}\frac1n\sum^{n-1}_{k=0}G\circ\Theta^k u=0, \,\,\mathbb{T} \,\,\,discrete,\]
\[\lim_{T\to\infty}\frac1T\int^{T}_{0}G\circ\Theta^t u\,dt=0, \,\,\mathbb{T} \,\,\,continuous,\]

for all $u\in B,$ and the convergence is uniform.
\end{lemma}

Applying this lemma to the continuous function $G=h-\rho,$ we have
$$\lim_{n\to\infty}\frac1{2n\pi}\sum^{n-1}_{k=0}h\circ\Theta^k(x,\omega)= \rho,$$
for all $(x,\omega)\in S^1\times\Omega.$ This completes the proof
of Proposition \ref{maprot3}
\end{proof}

\noindent{\bf Proof of Lemma \ref{mean}.} Since $B$ is a compact
metric space or a compact Hausdorff separable space, the space
$C(B)$ of the continuous functions on $B$ is separable. So we can
find a dense linear subspace $D$ generated by a countable set of
functions. If the statement is not true, then for some function
$G\in C(B)$, we can choose $D$ so that $G\in D$ and select
$n_i\in\mathbb{N}$(or $T_i\in\mathbb{R}$, respectively),\,$u_i\in
B$ such that $n_i(or T_i)\to\infty$ and
$$\lim_{n_i\to\infty}\frac1{n_i}\sum^{n_i-1}_{k=0}G\circ\Theta^k u_i=\delta\ne0,\mathbb{T} \,\,\,discrete,
$$
or respectively,
\[\lim_{T_i\to\infty}\frac1{T_i}\int^{T_i}_{0}G\circ\Theta^t u_i\,dt
=\delta\ne0, \,\,\mathbb{T} \,\,\,continuous.\]
 We may assume that $u_i\to\tilde u$. Using the Cantor diagonal process we can chose a
subsequence, which we still denote by $n_i(\text{or }T_i),\, u_i$,
such that
$$\lim_{n_i\to\infty}\frac1{n_i}\sum^{n_i-1}_{k=0}H\circ\Theta^k u_i\quad\text{ exists}$$
or respectively,
\[\lim_{T_i\to\infty}\frac1{T_i}\int^{T_i}_{0}H\circ\Theta^t u _i\,dt\quad\text{ exists}\]
for all $H\in D$. This limit defines a linear functional
$l=l(H),\,H\in D$.  Since $l$ is bounded with norm 1, it can be
extended uniquely to a bounded positive linear functional on
$C(B)$. Obviously, $l(H\circ\Theta^t)=l(H)$. By the Riesz
representation theorem, $l$ defines an invariant probability
measure $\mu$ on $B$. By our assumption
$$\int_B G\,d\mu=l(G)=\delta\ne0,$$
which is a contradiction. The proof of lemma is complete.

\hfill $\square$
\section{\bf Rotation Numbers of Random Differential Equations on the Circle}

In this section, we establish the analogous results of section 3
for random differential equations on the circle.\\

Let $(\Omega,\mathcal{F},\mathbb{P},(\theta^t)_{t\in\mathbb{R}})$
be a  metric dynamical system. Let $C(S^1)$ be the Banach space of
real-valued bounded continuous functions on
$S^1=\mathbb{R}/2\pi\mathbb{Z}$. Let $L^1(\Omega,C_L(S^1))$ denote
the set of the functions $f(\cdot,\omega)$ which is Lipschitz
continuous with respect to $x$ for each fixed $\omega$ and
$$M(\omega):=\max\left\{\sup_{x\in\mathbb{R}}|f(x,\omega)|,\sup_{x\ne
y}\frac{|f(x,\omega)-f(y,\omega)|}{|x-y|}\right\}\in
L^1(\Omega,\mathcal{F},\mathbb{P}).$$ In  $L^1(\Omega,C_L(S^1))$,
we introduce the norm
$$\|f\|=\int_\Omega\sup_x|f(x,\omega)|\,d\mathbb{P}.$$

Consider the following random ordinary differential equation on
the circle $S^1=\mathbb{R}/2\pi\mathbb{Z}$
\begin{equation}\label{equ1}
\dot x=f(x,\theta^t\omega).
\end{equation}
Let $x(t,x_0,\omega)$ denote the solution with the initial value
$x=x_0$ for $t=0.$  Assume that $f\in L^1(\Omega,C_L(S^1))$. Then,
by Theorem 2.2.6 in \cite{Arn98}, equation (\ref{equ1}) with
$x(0,x_0,\omega)=x_0$ has a unique global solution almost surely
in $\omega$.

\begin{theorem}\label{roteq} Let $f\in L^1(\Omega,C_L(S^1)).$
 Then,
\begin{itemize}
\item[(i)]there exists a $\theta^t$-invariant set
$\tilde\Omega\in\mathcal{F} $ of full measure and a $L^1$-function
$\rho(\omega):\tilde\Omega\to\mathbb{R}$ such that
$$\lim_{T\to\infty}\frac{x(T,x_0,\omega)-x_0}T=\lim_{T\to\infty}\frac1T\int^T_0f(x(t,x_0,\omega),\theta^t\omega)\,dt=\rho(\omega),\,\,
\text{ for all }x_0\in\mathbb{R},\,\,\omega\in\tilde\Omega,$$ and
$\rho(\theta^t\omega)=\rho(\omega)$ for all $t\in\mathbb{R},$
$\rho(\omega)$ is constant when $\theta^t$ is ergodic;

\item[(ii)] the map
$$\rho:   L^1_{\mathbb{P}}(\Omega,{\mathcal
L i p}_b(S^1))\to L^1(\Omega,\mathcal{F},
\mathbb{P}):f(x,\omega)\mapsto\rho(\omega)$$
 is continuous;

\item[(iii)] if, in addition, $\Omega$ is a compact metric space,
$\theta$  and $f(x,\omega)$ are continuous with $M(\omega)$
bounded, and $\mathbb{P}$ is the unique $\theta$-invariant
probability measure,
 then there exists $\rho\in \mathbb{R}$ such that
$$\lim_{T\to\infty}\frac{x(T,x_0,\omega)-x_0}T=\lim_{T\to\infty}\frac1T\int^T_0f(x(t,x_0,\omega),\theta^t\omega)\,dt=\rho\in\mathbb{R},\,\,$$
 for all $x_0\in\mathbb{R},\,\,\text{and all } \omega\in\Omega.$
\end{itemize}

\end{theorem}
\begin{proof}
We first prove (i). Since  $M(\omega)\in L^1$, we have that the
solution $x(t,x_0,\omega)$ of the initial value problem exists
globally for almost all $\omega\in\Omega$. We note that the
function $x(t,x_0,\omega)-x_0$ is $2\pi$-periodic
 and $x(t,x_0,\omega)$ is  increasing with respect to $x_0$.
 Therefore, the limit in (i), if exists, is independent of $x_0$. Let
$$\Omega_1:=\left\{\omega\in\Omega:\,\lim_{n\to\infty}\frac1n\int^{n+1}_nM(\theta^t\omega)dt=0\right\}.$$
From Lemma 2.2.5 and Proposition 4.1.3 in \cite{Arn98} it follows
that $\Omega_1$ is a $\theta^t$-invariant set of full measure. Let
$\tilde\Omega$ denote the set of $\omega\in\Omega_1$ such that the
limit in the (i) exists. Then, $\tilde\Omega$ is a
$\theta^t$-invariant set and the limit is $\theta$-invariant. Next,
we show that it has a full measure. By Proposition \ref{maprot1},
there exists an invariant set $\Omega_2$ of full measure such that
\[\lim_{\mathbb{N}\ni n\to\infty}\frac{x(n,x_0,\omega)-x_0}n
\] exists for $\omega\in\Omega_2$. For $\omega\in\Omega_1\cap
\Omega_2$, we have
$$\lim_{T\to\infty}\frac1T\int^T_0f(x(t,x_0,\omega),\theta_t\omega)\,dt=\lim_{T\to\infty}\frac1{[T]}\int^{[T]}_0f(x(t,x_0,\omega),\theta^t\omega)\,dt.$$
Since $\tilde\Omega\supset\Omega_1\cap\Omega_2$ and $\Omega_1$ and
$\Omega_2$ have a full measure, $\tilde \Omega$ has a full
measure.\\

Next, we show that (ii) holds. Denote by $x_f(t,x_0,\omega)$ the
solution of equation (\ref{equ1}) with the initial value $x=x_0$
at $t=0.$ Then, by Gronwall's  inequality, for $f,g \in
L^1_{\mathbb{P}}(\Omega,{\mathcal L i p}_b(S^1))$ we have

\[\|x_f(t,\cdot,\cdot)-x_g(t,\cdot,\cdot)\|\le t\|f-g\|\exp(t
(Lip(f))),\text{ for any given } t\in \mathbb{R}, \] which
together with Proposition (\ref{maprot2}) gives (ii) in the
theorem. \\

Finally, we show (iii). Let
$\Theta^t(x_0,\omega)=(x(t,x_0,\omega),\theta_t\omega)$,  the
corresponding skew product flow.  Since $\mathbb{P}$ is the unique
$\theta$-invariant probability measure, we have that $\theta$ is
ergodic and $\pi_\Omega\mu=\mathbb{P}$ for any invariant
probability  measure $\mu$ of $\Theta$. Therefore, by  (i) there
exists a constant $\rho$ such that for any invariant probability
measure $\mu$ of $\Theta$,
$$\rho=\lim_{T\to\infty}\frac1T\int^T_0f\circ\Theta^t(x_0,\omega)\,dt,\,\,a.s,\mu.$$
Thus,  by using Lebesgue's Dominated Convergence Theorem, Fubini's
Theorem, and Birkhoff's Ergodic Theorem, we have
\begin{align*}\rho&=\int\rho\,d\mu\\&=\int\left(\lim_{T\to\infty}\frac1T\int^T_0f\circ\Theta^t\,dt\right)\,d\mu\\&=
 \lim_{T\to\infty}\int\frac1T\int^T_0f\circ\Theta^t\,dt\,d\mu\\&= \lim_{T\to\infty}\frac1T\int^T_0\int f\circ\Theta^td\mu\,dt\,\\&=\int
 f\,d\mu.\end{align*}
Applying  Lemma \ref{mean} for the continuous time to the function
$f-\rho$, we obtain (3). This completes the proof of this theorem.
\end{proof}
\section{\bf Applications}

In this section, we apply our main results to almost periodic
differential equations. Let
$f(x,u):(\mathbb{R}/2\pi\mathbb{Z})\times\mathbb{R}^d\to
\mathbb{R}$ be a continuous function and $f(\cdot,u)$ be Lipschitz
continuous with respect to $x$ for each fixed $u$. Assume that
\[
M(u)=\lim_{\delta\to0}\sup_{0<|x-y|<\delta}\frac{|f(x,u)-f(y,u)|}{|x-y|}
\] is locally bounded measurable function in $u\in \mathbb{R}^d$. Suppose that $u(t):\mathbb{R}\to\mathbb{R}^d$ is an almost
periodic
 function. Consider the following almost periodic
time-depending differential equation on $S^1$
\begin{equation}\label{alm}
\dot x=f(x,u(t)).
\end{equation}
 Denote by $x(t,x_0)$ the solution of equation (\ref{alm}) with
 the initial value $x(0,x_0)=x_0$.

\begin{theorem}\label{app}
The rotation number $\rho$ defined by
 the limit
$$\rho:=\lim_{T\to\infty}\frac{x(T,x_0)-x_0}T=\lim_{T\to\infty}\frac1T\int^T_0f(x(t,x_0),u(t))\,dt$$
exists and is independent of $x_0$.
\end{theorem}
\begin{proof} Let $H(u):=\{u(\cdot+t):t\in\mathbb{R})$ be the
hull of $u$, i.e.,  the set of all translates of $u$. Then the
closure $\overline{H(f)}$ of $H(f)$ is compact and consists of
almost periodic functions, where the closure is taken  in the
uniform topology. $\overline{H(f)}$ has the structure of a compact
Abelian Polish group $G$ with unit $e=u$. The group operation $*$
is defined as follows: For $g=u(\cdot+t)$ and $h=u(\cdot+s)$,
$g*h=u(\cdot+t+s)$, while for $g=\lim u(\cdot+t_n)$ and $h=\lim
u(\cdot+s_n)$, $g*h=\lim u(\cdot+t_n+s_n).$

Associated to the almost periodic function $u$, we have  the
following canonical metric dynamical system. Let
$\mathbb{T}=\mathbb{R},\; \Omega:=G=\overline{H(f)},\;
\mathcal{F}$ the Borel $\sigma$-algebra of $G$. The metric
dynamical system is given by the translation of $\omega$ by $t$,
$\theta_t\omega:\omega(\cdot+t)$. Note that
$(t,\omega)\mapsto\theta_t\omega$ is continuous. The normalized
Haar measure $\mathbb{P}$ of $G$ is the unique $\theta$-invariant
probability. Under $\mathbb{P}$, $\theta$ is ergodic. We define a
continuous function $F$ on
$(\mathbb{R}/2\pi\mathbb{Z})\times\Omega$ as
$$F(x,\omega):=f(x,\omega(0)).$$ Now consider the random ordinary differential equation on the circle:
\begin{equation}\label{li3}\dot x=
F(x,\theta_t\omega).
\end{equation}
 Then by Theorem \ref{roteq}, for
each $\omega\in \Omega$, the rotation number of equation
(\ref{li3}) exists and is independent of $\omega\in\Omega.$
Obviously, for $\omega=u,\, F(x,\theta_t\omega)=f(x,u(t)).$ This
completes the proof.
\end{proof}

\vskip0.1in

\noindent {\bf Example 1.} Let
$A(x,u):\mathbb{R}^2\times\mathbb{R}^d\to\mathbb{R}^{2}$ be
continuous and be Lipschitz continuous in $x$ with a locally
bounded measurable Lipschitz constant in $u$ . Assume that $A$ is
positive homogeneous with respect to $x$, i.e.
$$A(\lambda x,u)=\lambda A(x,u),\,\,\,\text{ for }\lambda\ge0.$$
Let $u(t):\mathbb{R}\to\mathbb{R}^d$ be an almost periodic
function. Consider an almost periodic time-depending system in
$\mathbb{R}^2$ as
\begin{equation}\label{almsys}
\dot x=A(x,u(t)),\,\,\,x=(x_1,x_2)^T\in\mathbb{R}^2.
\end{equation}
In polar coordinates
$r=(x^2_1+x^2_2)^{1/2},\,\,\alpha=\arctan(x_2/x_1)$, the system
(\ref{almsys}) is written as
\[\dot r=\langle A(w,u(t)),w\rangle r,\quad \dot \alpha =\langle
A(w,u(t)),v\rangle,
\]
where
$w=(\cos\alpha,\sin\alpha)^T,\,\,v=(-\sin\alpha,\cos\alpha)^T.$
Denote by $(r(t,r_0,\alpha_0),\alpha(t,\alpha_0))$ the solution of
the system with the initial value condition
$(r(0,r_0,\alpha_0),\alpha(0,\alpha_0))=(r_0,\alpha_0)$. Then the
rotation number of the system (\ref{almsys}) is defined to be the
linear growth rate of the angular component, i.e.,  by
\[\rho:=\lim_{T\to\infty}\frac{\alpha(t,\alpha_0)-\alpha_0}T=\lim_{T\to\infty}\frac1T\int^T_0\langle
A(w,u(t)),v\rangle\,dt.\] By the Theorem \ref{app}, the rotation
number exists and is independent of the initial value.

\vskip0.1in

\noindent {\bf Example 2.} Let $f_1,f_2,...,f_k$ be the
orientation preserving homeomorphisms of the circle with the same
invariant probability measure $\mu_x$ and let $(p_1,p_2,...,p_k)$
be a probability vector with non-zero entries (i.e., $p_i>0$ for
each $i$ and $\sum^k_{i=1}p_i=1$). Assume that $f_i$ has the
rotation number $\rho_i$. Let $(K,2^K,\mu)$ denote the probability
space where $K=\{1,2,...,k\}$ and the point $i$ has measure $p_i.$
Let \[(\Omega, \mathcal{F},
\mathbb{P})=\prod^\infty_1(K,2^K,\mu).\] We write points of
$\Omega$ in the form $\omega=(\omega_1,\omega_2,...),\omega_i\in
K,$  and define $\theta:\Omega\to\Omega$ by
\[\theta(\omega_1,\omega_2,...)=(\omega_2,\omega_3,...).\] Then,  $\theta$
is $\mathbb{P}$ measure-preserving. We define an orientation
preserving random map of the circle as
$\varphi(x,\omega)=f_{\omega_1}.$ Then the product measure
$\mu_x\times\mathbb{P}$ is the invariant measure of dynamical
system $\Theta^n(x, \omega)=(\phi(n,\omega)x, \theta^n\omega))$,
where
\[\phi(n,\omega)x=f_{\omega_n}\circ \cdots \circ f_{\omega_2} \circ f_{\omega_1}(x)\] is the random
dynamical
system over $\theta$
generated by $\varphi(x,\omega)$.

Let $f_i(x)=x+h_i(x), \varphi(x,\omega)=x+h(x,\omega)$. By formula
(6), the rotation number of $\phi$ is
\begin{align*}\rho&=\frac1{2\pi}\int_{S^1\times\Omega}hd\mathbb{P}d\mu_x\\
&=\frac1{2\pi}\int_{S^1}\sum^k_{i=1}p_ih_i(x)d\mu_x\\
&=\frac1{2\pi}\sum^k_{i=1}p_i\lim_{n\to\infty}\int_{S^1}\frac1n\sum^{n-1}_{m=0}h_i\circ
f^m_id\mu_x\\
&=\frac1{2\pi}\sum^k_{i=1}p_i\int_{S^1}\lim_{n\to\infty}\frac1n\sum^{n-1}_{m=0}h_i\circ
f^m_id\mu_x\\
&=\frac{1}{2\pi}\sum^k_{i=1}p_i\int_{S^1}2\pi\rho_id\mu_x\\
&=\sum^k_{i=1}p_i\rho_i.
\end{align*} This example shows that the rotation number of random
compositions of the orientation preserving homeomorphisms of the
circle, $f_1, \cdots, f_k$, is the probability weighted average of
the rotation numbers of them.

\section{\bf Analytical Conjugacy to a Circle Rotation}

In this section, we study the problem of analytic conjugacy to a
circle rotation. We first review a Diophantine condition.

\begin{definition} We say that $\mu\in \mathbb{R}^m$ is a vector of
type $(C,\nu)$ if
\begin{equation}\label{cnu}|e^{2\pi i\langle\mu,\,k\rangle}-1|>\frac
C{|k|^\nu}, \quad |k|:=|k_1|+|k_2|+\cdots+|k_m|
\end{equation}
for all nonzero integer vector $k\in\mathbb{Z}^m.$
\end{definition}

The next lemma gives that almost all of vectors $\mu\in
\mathbb{R}^m$ satisfy (\ref{cnu}).
\begin{lemma} Let $\nu>m$ be a constant. For almost every real
vector $\mu\in \mathbb{R}^m$, there exists  $C=C(\mu,\nu)>0$ such
that the inequality (\ref{cnu}) holds for all nonzero integer
vector $k\in \mathbb{Z}^m$.
\end{lemma}
\begin{proof} The proof of this lemma follows from the standard argument.
First we claim that for almost every vector
$\mu\in\mathbb{R}^m$ there exists $C=C(\mu,\nu)>0$ such that
\begin{equation}\label{neq}|\langle k,\mu\rangle-q|>\frac C{|k|^\nu}
\end{equation}
 for all
$k\in\mathbb{Z}^m\backslash\{0\}$ and $q\in\mathbb{Z}$. Indeed, we
fix a ball in $\mathbb{R}^m$ and estimate the measure of the set
of $\mu$ in it which does not satisfy the inequality (\ref{neq}).
Let
\[L_{k,q}=\{\mu\in\mathbb{R}^m:\,\langle k,\mu\rangle-q=0\}\]
denote the resonance plane. The inequality
\[
|\langle k,\mu\rangle-q|\leq \frac C{|k|^\nu}
\]determines a neighborhood of width not greater than
$C_1C/|k|^{\nu+1}$ of the resonance plane. Therefore, the measure
of the part of this neighborhood which is contained in the ball
does not exceed $C_2C/|k|^{\nu+1}$. Summing over $k$ with fixed
$|k|$, we obtain that the measure is not more than
$C_3C/|k|^{\nu-m+2}.$ Summing over $q$ with fixed $|k|$ such that
distance between $L_{k,q}$ and the ball is less than $1$, we
obtain the measure is bounded by  $C_4C/|k|^{\nu-m+1}$. Summing
over $|k|$, we have that the measure is bounded by
$C_5(\nu)C<\infty.$ Consequently, the set of $\mu$ in the ball is
covered by the sets of arbitrarily small measure. Hence, such a
set has measure zero.

The proof of this lemma  follows this claim since the distance of
$\langle k,\mu\rangle$ from the closest integer is bounded from
below by $C/|k|^\nu$ and a chord of the unit circle is not shorter
than the length of the small arc subtended by it divided by
$\frac\pi2$. This completes the proof of the lemma.
\end{proof}
\vskip0.1in

We now consider a class of random maps of the circle over an $m-1$
dimensional torus. Let $\Omega=\mathbb{R}^{m-1}/2\pi\mathbb{Z}$ be
the torus of dimension $m-1$. Consider the metric dynamical system
$\{\theta^{n}\}_{n\in\mathbb{Z}}$ on $\Omega$ given by
\[\theta^n:\omega\mapsto\omega+2n\pi\alpha,\]
where $\alpha\in\mathbb{R}^{m-1}$ is a given vector. We assume
that
$$\langle\alpha,k\rangle-j\ne0,\,\,\text{for all }k\in\mathbb{Z}^{m-1}\backslash\{0\},\,j\in\mathbb{Z}.$$
Then,  the normalized Lebesgue measure $\mathbb{P}$ is the unique
$\theta$-invariant probability measure and $\theta$ is ergodic
under $\mathbb{P}$. Let
$\varphi(x,\omega):(\mathbb{R}/2\pi\mathbb{Z}\times\Omega)\to
\mathbb{R}/2\pi\mathbb{Z}$ be an orientation preserving random map
of the circle over $\theta$. Suppose that $\varphi(\cdot,\cdot)$
is continuous. Then,  by Theorem A, the rotation number of
$\varphi$ exists and is given by
\[\rho=\lim_{n\to\infty}\frac{\phi(n, \omega)x}{2\pi n}=\lim_{n\to\infty}\frac1{2\pi n}
\varphi(\cdot, {\theta^{n-1}\omega})\circ \varphi(\cdot,
{\theta^{n-2} \omega})\circ\cdots\circ \varphi(x,\omega), \quad
\text{for all } x\in \mathbb{R}, \omega\in \Omega, \] which is
independent of $\omega$
and $x$. \\


Consider a perturbation of the circle rotation by $2\pi \rho$:
\[
\varphi(x, \omega)=x+2\pi \rho +p(x, \omega),
\] where $p(x,\omega)$ is a holomorphic function defined on the strip
$U_r$ which was introduced in the introduction.\\

We have the following theorem on analytic conjugacy to a circle
rotation.

\vskip0.1in

\begin{theorem}\label{thmc1} { Let $p(x,\omega)$ be analytic in $U_r$ and
$2\pi$-periodic in each variable, real on the real axes. Assume

\begin{itemize}

\item[(1)] $\varphi(x,\omega)=x+2\pi\rho+p(x,\omega)$ has the
rotation number $\rho$ and

\item[(2)] the vector $\mu=(\rho,\alpha)$ is of $(C,\,\nu)$ type,
i.e.,
\begin{equation}\label{cnu}|e^{2\pi i\langle\mu,\,k\rangle}-1|>\frac
C{|k|^\nu}, \quad |k|:=|k_1|+|k_2|+\cdots+|k_m|
\end{equation}
for all nonzero integer vector $k\in\mathbb{Z}^m,$ where $c$ and
$\nu$ are positive constants.
\end{itemize}
Then, there exists $\epsilon>0$ depending only on $C,\nu, r$ and
$m$ such that if $\|p\|_r<\epsilon$, then the random map
$\varphi(\cdot, \omega)$ is analytically conjugate to the circle
rotation by the angle $2\pi\rho$, i.e.,  there exists an
analytical random transformation
\[H(\cdot,\cdot):\mathbb{R}^m/2\pi\mathbb{Z}\to\mathbb{R}/2\pi\mathbb{Z}\]
such that
\[H(x+2\pi\rho, \theta\omega)=\varphi (\cdot, \omega)\circ H(x, \omega).\]
}
\end{theorem}
\vskip0.1in

Instead of proving Theorem \ref{thmc1}, we will prove a more
general result which gives Theorem \ref{thmc1} as its special
case. We consider a mapping
\[
\Phi(z):\mathbb{R}^m/2\pi\mathbb{Z}\to\mathbb{R}^m/2\pi\mathbb{Z}.
\]
\begin{definition} We say that the mapping $\Phi(z)$  has a rotation vector $\mu\in\mathbb{R}^m$, if
\[\lim_{n\to\infty}\frac1n\Phi^{n}(z)=2\pi\mu,\]
for all $z\in \mathbb{R}^m/2\pi\mathbb{Z}$.
\end{definition}

Consider a perturbation of a vector rotation:
\[\Phi(z)=z+2\pi\mu+p(z):\mathbb{R}^m/2\pi\mathbb{Z}\to\mathbb{R}^m/2\pi\mathbb{Z}\]
We have the following result.
\begin{theorem}\label{conj2} Let $p(z)$ be
analytic in $U_r$ and $2\pi$-periodic function in each variable,
real on the real axes. Assume
\begin{itemize}
\item[(1)]
 $\Phi(z)=z+2\pi\mu+p(z):\mathbb{R}^m/2\pi\mathbb{Z}\to\mathbb{R}^m/2\pi\mathbb{Z}$
has a rotation vector $\mu$ and

\item[(2)] the rotation vector $\mu$ is of $(C,\nu)$ type for
positive constants $C$ and $\nu$.
\end{itemize}
Then, there exists $\epsilon>0$ depending only on $C,\nu, r$ and
$m$ such that if $\|p\|_r<\epsilon$, then $\Phi(z)$ is
analytically conjugate to the rotation $R_\mu:z\mapsto z+2\pi\mu$,
i.e., there exists  an analytic transformation
\[
H(z)=z+h(z):\mathbb{R}^m/2\pi\mathbb{Z}\to\mathbb{R}^m/2\pi\mathbb{Z},
\]
where $h(z)$ is a $2\pi$-periodic function in each variable, such
that
\begin{equation}\label{conjeq}
H\circ R_\mu=\Phi\circ H.
\end{equation}
\end{theorem}
\vskip0.1in

The proof of this theorem is based on the classic KAM approach. We
write $H(z)$ in the form $H(z)=z+h(z)$, where $h(z)$ is
$2\pi$-periodic in each variable.  Then substituting it into
conjugate equation (\ref{conjeq}), we obtain the functional
equation for $h$
\[h(z+2\pi\mu)-h(z)=p(z+h(z)).\]
The first approximation of this equation is the so called {\it
homological equation} for $h$,
\begin{equation}\label{homoeq}
h(z+2\pi\mu)-h(z)=p(z).
\end{equation}
Obviously, the homological equation is not solvable if the mean
value of $p(z)$ is nonzero. We expand the given function $p$ and
the unknown function $h$ in the Fourier series:
\[p(z)=\sum_{k\in\mathbb{Z}^m\setminus\{0\}}p_ke^{i\langle k,\,z\rangle},\qquad
h(z)=\sum_{k\in\mathbb{Z}^m\setminus\{0\}}h_ke^{i\langle
k,\,z\rangle}.\] Plugging them into (\ref{homoeq}) and comparing
the coefficients of $e^{i\langle k,\,z\rangle}$, we have
\begin{equation} \label{coef}
h_k=\frac{p_k}{e^{2\pi i\langle \mu,\,k\rangle}-1},\quad
k\in\mathbb{Z}^m\setminus\{0\}.
\end{equation}

The idea of the proof of Theorem \ref{conj2} is the following. We
solve the homological equation (\ref{homoeq}) with the right-hand
side $\tilde p(z)=p(z)-p_0$, where $p_0$ is the mean value of the
function $p(z)$. Denote by $h^0$ the solution. Let
$H_0(z)=z+h^0(z)$. Set $\Phi_1=H^{-1}_0\circ \Phi\circ H_0$ and
define a function $p^1(z)$ by the relation
$$\Phi_1(z)=z+2\pi\mu+p^1(z).$$
The next approximation is constructed in the same way.  Beginning
with $\Phi_1$ in the place of $\Phi$, we solve the corresponding
homological equation for $h^1$ and  let $H_1=z + h^1(z)$. The
transformation  $H_1$ converts $\Phi_1$ into
$$\Phi_2=H^{-1}_1\circ \Phi_1\circ H_1.$$
Repeating this procedure,  we have a sequence of transformation
$H_n.$ Let
$$\textbf{H}_n=H_0\circ H_1\circ\cdots\circ H_{n-1}.$$
Then
$$\Phi_n=\textbf{H}^{-1}_n\circ \Phi\circ\textbf{H}_n. $$
Finally, we will prove that $\lim_{n\to\infty}\Phi_n=z+2\pi\mu.$\\

Before proving the theorem, we first introduce several lemmas
which we need later.
\begin{lemma}\label{le1} Let $f: \mathbb{R}^m \to\mathbb{R}$
be $2\pi$-periodic in each component, and   be analytic in the
strip $U_r$ and continuous in the closure of this strip. Assume
that $\|f\|_r\le M$. Then its Fourier coefficients satisfy
\[|f_k|\le Me^{-|k|r}.\]
\end{lemma}
\begin{proof} For $k=(k_1,k_2,\dots,k_m)\in\mathbb{Z}^m$, we
define $u=(u_1,u_2,\dots,u_m)\in\mathbb{Z}^m$ as
$u_j=-\text{sgn}(k_j),\,j=1,2,\dots,m.$ Then we have
\[f_k=\frac1{(2\pi)^m}\int_{\mathbb{T}^m}e^{-i\langle k,z\rangle} f(z)\,dz=\frac1{(2\pi)^m}\int_
{\mathbb{T}^m}e^{-i\langle k,z\rangle-|k|r}f(z+iru)\,dz.\] Hence,
\[|f_k|\le\frac1{(2\pi)^m}\int_
{\mathbb{T}^m}|e^{-i\langle k,z\rangle-|k|r}f(z+iru)|\,dz\le
Me^{-|k|r}.\]
\end{proof}
\begin{lemma}\label{le2} If $|f_k|\le Me^{-|k|r}$, then the
function $f=\sum f_ke^{i\langle k,\,z\rangle}$ is analytic in the
strip $U_r$ and $\|f\|_{r-\delta}\le
8\left(\frac{4m-4}e\right)^{m-1}M\delta^{-m},$ where
$\delta<\min\{1,r\}.$
\end{lemma}
\begin{proof} The proof follows from the following computation.
\begin{align*}
\|f\|_{r-\delta}&\le\sum|f_k||e^{i\langle k,\,z
\rangle}|\le \sum Me^{-|k|r}e^{|k|(r-\delta)}\\&=M\sum
e^{-|k|\delta}=M\sum^{\infty}_{l=0}\frac{2^m(l+m-1)!}{l!(m-1)!}e^{-l\delta}
\\&\le2^mM\sum^{\infty}_{l=0}(l+1)^{m-1}e^{-l\delta}\le2^mMe^\delta\left(\frac{2m-2}{e\delta}\right)^{m-1}
\sum^{\infty}_{l=0}e^{-\delta l/2}\\&\le
Me^\delta\left(\frac{4m-4}{e}\right)^{m-1}\frac1{\delta^{m-1}}\frac1{1-e^{-\delta/2}}\le8M\left(\frac{4m-4}e\right)^{m-1}
\delta^{-m}.
\end{align*}
\end{proof}

\begin{lemma}\label{le3}Let  $p(z):U_r\to\mathbb{C}^m$ be a $2\pi$-periodic
analytic function  with mean value $0$. Let $h(z)$ be the solution
of homological equation (\ref{homoeq}). Then, there exists a
constant $\lambda=\lambda(C,\nu,m)>0$ such that if $\mu$ is of
type $(C,\nu)$, then for any $r<\frac12$ and  any $\delta>0$
smaller than $r$, we have
$\|h\|_{r-\delta}\le\|a\|_r\delta^{-\lambda}.$
\end{lemma}
\begin{proof} Let $M=\|p\|_r$,  $p(z)=(p_1(z),p_2(z),\dots,p_m(z)),$ and $ h(z)=(h_1(z),h_2(z),\dots,h_m(z))$.
We write
\[p_j(z)=\sum_{k\in\mathbb{Z}^m\setminus\{0\}}p^j_ke^{i\langle k,\,z\rangle},\qquad
h_j(z)=\sum_{k\in\mathbb{Z}^m\setminus\{0\}}h^j_ke^{i\langle
k,\,z\rangle},\,\, j=1,2,\dots,m.\] By Lemma \ref{le1}, we have
$|p^j_k|\le Me^{-|k|r}.$ Since $\mu$ is of type $(C,\nu)$, using
(\ref{coef}) we have
$$|h^j_k|\le |k|^\nu Me^{-|k|r}/C\le MC^{-1}|k|^\nu
e^{-\delta|k|/2}e^{-|k|(r-\delta/2)}\le MC^{-1}\left(\frac\nu
e\right)^{\nu}\left(\frac\delta2\right)^{-\nu}e^{-|k|(r-\delta/2)}.$$
By Lemma \ref{le2},
\[\|h\|_{r-\delta}=\max_{1\le j\le m}\|h_j\|_{r-\delta}\le8MC^{-1}\left(\frac\nu
e\right)^{\nu}\left(\frac\delta2\right)^{-\nu}\left(\frac{4m-4}e\right)^{m-1}\left(\frac\delta2\right)^{-m}\le
M\delta^{-\lambda},\] for $\lambda$ sufficiently large.
\end{proof}
The next two lemmas are obvious.
\begin{lemma}\label{le4} Let $H(z)=z+h(z):\mathbb{R}^m/2\pi\mathbb{Z}\to\mathbb{R}^m/2\pi\mathbb{Z}$
be a homeomorphism, where $h$ is $2\pi$-periodic in each component
of $z$. Suppose that the mapping
$\Phi:\mathbb{R}^m/2\pi\mathbb{Z}\to\mathbb{R}^m/2\pi\mathbb{Z}$
has rotation vector $\mu\in\mathbb{R}^m$, then the mapping
$H^{-1}\circ \Phi\circ H$ has also the rotation vector $\mu$.
\end{lemma}
\begin{lemma}\label{le5} Suppose that the homeomorphism
\[\Phi(z)=z+2\pi\mu+p(z):\mathbb{R}^m/2\pi\mathbb{Z}\to\mathbb{R}^m/2\pi\mathbb{Z}\]
has the rotation vector $\mu$, where $p$ is $2\pi$-periodic in
each component of $z$.  Then each component $p_j(z)$ of $p(z)$
vanishes at some point.
\end{lemma}
The next one is our main lemma for proving Theorem \ref{conj2}.
\begin{lemma}\label{le6} There exist positive constants $\tau,\gamma$
depending only on $C,\nu$ and $m$ such that for every $\delta$ in
the interval $(0,r)$, where $r<1/2$, we have
\[\|p^1\|_{r-\delta}\le\|p\|^2_r\delta^{-\gamma},\,\text{ provided that } \|p\|_r\le\delta^\tau.\]
\end{lemma}
\begin{proof} We first show that if $\tau$ is large enough, then
$A_1$ is well defined and analytic in the strip $U_{r-\delta}$.

Let $M=\|p\|_r$ and assume $M\le\delta^\tau$. Then, the mean $p_0$
of $p$ satisfies $\|p_0\|\le M$ and $\|\tilde p\|_r=\|p-p_0\|_r\le
2M.$ By Lemma \ref{le3}, for any $0<\alpha<r$, we have
$\|h^0\|_{r-\alpha}\le2M\alpha^{-\lambda}$, which implies
\begin{equation}\label{ineq-1}
\|Dh^0\|_{r-2\alpha}\le2mM\alpha^{-\lambda-1}.
\end{equation}
 Let
$\alpha=\delta/8$. Choosing $\tau$ is sufficiently large, we
obtain
\[\|p\|_r<\alpha,\quad\|h^0\|_{r-\alpha}<\alpha,\quad\|Dh^0\|_{r-2\alpha}<\alpha.\]
Therefore, $H_0(z)=z+h(z)$ is a diffeomorphism on $ U_{r-2\alpha}$
and its image contains $U_{r-3\alpha}$. Since
$H_0U_{r-\delta}\subset U_{r-\delta+\alpha}$, $\Phi\circ
H_0U_{r-\delta}\subset U_{r-\delta+2\alpha}\subset U_{r-3\alpha}.$
The inverse $H^{-1}_0$ is defined on $\Phi\circ H_0U_{r-\delta}.$
Hence, the mapping $\Phi_1=H^{-1}_0\circ \Phi\circ H_0$ is well
defined and analytic on $U_{r-\delta}$

Next,  we estimate the function $p^1$. Since $H_0\circ
\Phi_1=\Phi\circ H_0,$ we have
\[z+2\pi\mu+p^1(z)+h^0(z+2\pi\mu+p^1(z))=z+h^0(z)+2\pi\mu+p(z+h^0(z)),\]
or
\begin{align*}p^1(z)&=p(z+h^0(z))+h^0(z)-h^0(z+2\pi\mu+p^1(z))\\
&=p(z+h^0(z))-p(z)+h^0(z+2\pi\mu)-h^0(z+2\pi\mu+p^1(z))+p_0.
\end{align*}
Thus,
\begin{align} \label{ineq0}
\begin{split}
\|p^1(z)\|_{r-\delta}&\le
\|p(z+h^0(z))-p(z)\|_{r-\delta}+\|h^0(z+2\pi\mu)-h^0(z+2\pi\mu+p^1(z))\|_{r-\delta}+\|p_0\|
\end{split}
\end{align}
 By the mean value theorem and the Cauchy inequality, we have
\begin{equation}\label{ineq1}
\|p(z+h^0(z))-p(z)\|_{r-\delta}\le\|Dp\|_{r-\alpha}\|h^0\|_{r-\delta}\le\frac
M\alpha M\delta^{-\lambda}\\ =8M^2\delta^{-\lambda-1}\le
M^2\delta^{-\kappa},
\end{equation}
where the constant $\kappa$ depends only on $\lambda$, i.e., hence
only on $C,\nu$ and $m$. Similarly, we have
\begin{equation}\label{ineq2}
\|h^0(z+2\pi\mu)-h^0(z+2\pi\mu+p^1(z))\|_{r-\delta}\le\|Dh^0\|_{r-\alpha}\|p^1\|_{r-\delta}
\le2mM\alpha^{-\lambda-1}\|p^1\|_{r-\delta}. \end{equation}
By Lemma \ref{le4}, the rotation vector of $\Phi_1$ is $\mu$.
Thus, from Lemma \ref{le5}, every component $p^1_j(z)$ of $p^1(z)$
vanishes at some real point, say $z^j_0\in\mathbb{R}^m,\,
j=1,2,\dots,m$. Let $p_0=(p^0_1,p^0_2,\dots,p^0_m).$ Then
\[p^0_j=p_j(z^j_0)-p_j(z^j_0+h^0(z^j_0))+h^0_j(z^j_0+2\pi\mu+p^1(z^j_0))-h^0_j(z^j_0+2\pi\mu).\]
Hence,
\begin{align} \label{ineq3}
\begin{split}\|p_0\|&=\max_{1\le j
\le m}|p^0_j|\\&\le\max_{1\le j\le
m}\{|p_j(z^j_0)-p_j(z^j_0+h^0(z^j_0))|+|h^0_j(z^j_0+2\pi\mu+p^1(z^j_0))-h^0_j(z^j_0+2\pi\mu)|\}\\&
\le\|p(z)-p(z+h^0(z))\|_{r-\delta}+\|h^0(z+2\pi\mu+p^1(z))-h^0(z+2\pi\mu)\|_{r-\delta}
\\&\le
M^2\delta^{-u}+2mM\alpha^{-\lambda-1}\|p^1\|_{r-\delta}.
\end{split}
\end{align}
Combining \ref{ineq1}),(\ref{ineq2}), and (\ref{ineq3}) with
(\ref{ineq0}), we obtain
\begin{align*}\|p^1\|_{r-\delta}&\le
2M^2\delta^{-\kappa}\left(1-4mM\alpha^{-\lambda-1}\right)^{-1}\le2M^2\delta^{-\kappa}\left(1-4m\delta^\tau(\delta/8)^{-\lambda-1}\right)^{-1}\\
&=2M^2\delta^{-\kappa}\left(1-m2^{3\lambda+5}\delta^{\tau-\lambda-1}\right)^{-1}\le4M^2\delta^{-\kappa}\le
M^2\delta^{-\gamma},
\end{align*}
provided that $\tau$ is sufficiently large. Here
$\gamma=\kappa+2$.
\end{proof}
{\it Proof of Theorem \ref{conj2}.} Let
$\delta_0\le\frac12,\,\delta_n=\delta^{3/2}_{n-1}$, and fix
$\delta_0$ small enough such that
$\sum^{\infty}_{n=0}\delta_n<\frac r2$. Set
$r_0=r,\,r_n=r_{n-1}-\delta_{n-1},\, M_n=\delta^N_n$, where
$N=\max\{\tau,2\lambda\}.$ Assume that $\|p\|_r\le M_0$, we claim
that
\begin{equation}\label{ineq4}
\|p^n\|_{r_n}\le M_n,\quad n\ge 0.
\end{equation}
We prove (\ref{ineq4}) by induction. Assume that for $n=k$,
inequality (\ref{ineq4}) holds. Then by Lemma \ref{le6},
\[\|p^{k+1}\|_{r_{k+1}}\le M^2_k\delta^{-\gamma}_k=\delta^{2N-\gamma}_k\le\delta^{3N/2}_k=\delta^N_{k+1}=M_{k+1}.\]
Next,  we prove the convergence of the composition
$\textbf{H}_n=H_0\circ H_1\circ\cdots\circ H_{n-1}$. First, by
Lemma \ref{le3} and (\ref{ineq-1}), we have that for $\tau$
sufficiently large the diffeomorphism $H_0$ is analytic in $U_r$
and satisfies $\|h^0\|_{r_1}<\delta_0,\,\|Dh^0\|_{r_1}<\delta_0.$
By induction, we have that for any $n\ge0,$
$$\|h^{n-1}\|_{r_n}<\delta_{n-1},\|Dh^{n-1}\|_{r_n}<\delta_{n-1}$$
and $\textbf{H}_n$ is analytic in $U_{r_n}$. Therefore,  the
derivative of $\textbf{H}_n$ satisfy
\[0<C_-:=\prod^\infty_{k=0}(1-\delta_k)\le\|D\textbf{H}_n\|_{r_n}\le\prod^\infty_{k=0}(1+\delta_k):=C_+.\]
Hence $\textbf{H}_n$ is a diffeomorphism. Note that
\[\|\textbf{H}_n-\textbf{H}_{n+1}\|_{r/2}\le C_+\|h^n\|_{r/2}\le C_+\delta_n.\]
So,   the sequence $\textbf{H}_n$ converges on $U_{r/2}.$ Let
\[\textbf{H}=\left.\lim_{n\to\infty}\textbf{H}_n\right|_{U_{r/2}}.\]
Then,
\[\textbf{H}\circ R_\mu=\lim_{n\to\infty}\textbf{H}_n\circ \Phi_n=\lim_{n\to\infty}\Phi\circ\textbf{H}_n
=\Phi\circ \textbf{H}.\]
This completes the proof of this theorem.

\hfill$\square$


\begin{thebibliography}{1}

\bibitem{Arn98}
L.~Arnold.
\newblock {\em {R}andom {D}ynamical {S}ystems}.
\newblock Springer, New York, 1998.

\bibitem{Arnold1}
V. I. ~Arnold.
\newblock {\em Geometric Methods in the Theory of Ordinary
Differential. Equations,}
\newblock Springer-Verlag, New York, 1983.


\bibitem{Arnold2} V. I. Arnold.
\newblock Small denominators I. One the mapping of a circle into
itself.
\newblock{\em Izv. Akad. Nauk. Math.}, 25 (1961), 21-86.




\bibitem{CLS}  S-N. ~Chow, K. ~Lu, and Y-Q. ~Shen.
\newblock Normal form and linearization for quasiperiodic systems.
\newblock {\em Trans. Amer. Math. Soc.} 331 (1992), 361--376.



\bibitem{Johnson1} R. Fabbri, R. Johnson, and C.  Núñez.
\newblock Rotation number for non-autonomous linear Hamiltonian systems.
II. The Floquet coefficients, \newblock{\em  Z. Angew. Math.
Phys.}, 54 (2003),  652--676.

\bibitem{Johnson2} R. Fabbri, R. Johnson, and C.  Núñez.
\newblock Rotation number
for non-autonomous linear Hamiltonian systems. I. Basic
properties. \newblock {\em Z. Angew. Math. Phys.}, 54 (2003),
484--502.

\bibitem{Johnson3} R. Fabbri, R. Johnson, and C.  Núñez.
\newblock Disconjugacy and the rotation number for linear, non-autonomous
Hamiltonian systems, \newblock {\em Ann. Mat. Pura Appl.}, 185
(2006),  S3--S21.

\bibitem{Zhang} H. Feng and M.  Zhang. \newblock Optimal estimates on rotation number of
almost periodic systems, \newblock{\em  Z. Angew. Math. Phys.}, 57
(2006), 183--204.



\bibitem{Franks} J. Franks. \newblock Rotation numbers and instability
sets,
\newblock {\em  Bull. Amer. Math. Soc., }  40 (2003),  263--279.


\bibitem{Herman} M. R. Herman.
\newblock Sur la conjugaison differentielle des diff\'eomorphismes
du cercle \'a des rotations,
\newblock {\em Publ. Math. I. H. E. S.,} 49(1979), 5-234.

\bibitem{JohnsonMoser} R.  Johnson and J. Moser. \newblock The rotation number for
almost periodic potentials, \newblock{\em Comm. Math. Phys.},  84
(1982),  403--438.

\bibitem{KatokHass} A. Katok and B.  Hasselblatt. \newblock Introduction to the modern
theory of dynamical systems. With a supplementary chapter by Katok
and Leonardo Mendoza. Encyclopedia of Mathematics and its
Applications, \newblock {\em 54. Cambridge University Press,
Cambridge}  1995.

\bibitem{KO} Y. Katznelson and D. Ornstein.
\newblock The differentiability of the conjugation of certain
diffeomorphism of the circle,
\newblock{\em Ergod. Th. and Dyn. Sys.,} 9 (1989), 643-680.




\bibitem{Moser} J. K. ~Moser.
 \newblock  A rapidly covergent iteration method and nonlinear differential equations
 II.
 \newblock {\em Ann. Scuo. Norm. Sup. Pisa.} 20 (1966), 499--535



\bibitem{Poincare} H. ~Poincar\'e.
 \newblock Sur le probl\'eme des trois corps et les \'equations de la dynamique
 \newblock {Acta Math.}, 13 (1890), 1--270.

\bibitem{Ruffino} P.R.C. Ruffino.
\newblock A sampling theorem for rotation numbers of linear
processes in $\mathbb{R}^2$, Random operators and Stochastic
Equations, 8(2000), 175-188.


\bibitem{SK} Ya. G. Sinai and K. M. Khanin.
\newblock Smoothness of conjugacies of diffeomorphisms of the
circle with rotations,
\newblock{\em Russ. Math. Surv.,} 44 (1989), 69-99.




\bibitem{Yoccoz1} J.-C. Yoccoz.
\newblock Conjugaison differentielle des diff\'eomorphismes du
cercle dont le nombre de rotation v\'erifie une condition
Diophantienne,
\newblock{\em Ann. Sci. Ec. Norm. Sup.}, 17 (1984), 333-361.





\end{thebibliography}
\end{document}